\documentclass[MSNbibl,number,citesort,secthm,seceqn,dvips]{arxbj}
\usepackage{upgreek}

\aid{0}
\volume{18}
\issue{4}
\pubyear{2012}
\firstpage{1172}
\lastpage{1187}
\doi{10.3150/11-BEJ382}

\makeatletter
 \newtheorem{lem}[thm]{Lemma}
\newtheorem{cor}[thm]{Corollary}
 \newtheorem{prop}[thm]{Proposition}
\newremark{ex}[thm]{Example}
 \newremark{rem}[thm]{Remark}

\newcommand{\cL}{\mathcal{L}}
\newcommand{\al}{\alpha}
\newcommand{\Gm}{\Gamma}
\newcommand{\ld}{\lambda}
\newcommand{\wt}{\widetilde}
\newcommand{\R}{\mathbb{R}}
\newcommand{\RR}{\mathbb{R}}
\newcommand{\N}{\mathbb{N}}
\newcommand{\LL}{\mathbb{L}}
\newcommand{\law}{\mathcal L}
\newcommand{\overset}{\stackrel}
\newcommand{\eqref}[1]{(\ref{#1})}
\renewcommand{\epsilon}{\varepsilon}
\newcommand{\fracd}[2]{({#1}/{#2})}
\newcommand{\fraca}[2]{{#1}/{#2}}
\newcommand{\fracc}[2]{{#1}/{(#2)}}
\makeatother

\begin{document}
\begin{frontmatter}

\title{Distributions of
exponential integrals of independent increment processes related to generalized gamma convolutions}
\runtitle{Exponential integrals related to GGCs}

\begin{aug}
\author[1]{\fnms{Anita} \snm{Behme}\corref{}\thanksref{1}\ead[label=e1]{a.behme@tu-bs.de}},
\author[2]{\fnms{Makoto} \snm{Maejima}\thanksref{2}},
\author[3]{\fnms{Muneya}~\snm{Matsui}\thanksref{3}} \and
\author[4]{\fnms{Noriyoshi}~\snm{Sakuma}\thanksref{4}}

\runauthor{Behme, Maejima, Matsui and Sakuma} 

\address[1]{TU Braunschweig, Pockelsstr. 14,
38106 Braunschweig, Germany. \printead{e1}}
\address[2]{Keio University, 3-14-1, Hiyoshi, Kohoku-ku, Yokohama
223-8522, Japan}
\address[3]{Department of Business Administration, Nanzan University,
18 Yamazato-cho, Showa-ku, Nagoya 466-8673, Japan}
\address[4]{Department of Mathematics Education, Aichi
University of Education, 1 Hirosawa, Igaya-cho, Kariya-shi, 448-8542, Japan}
\end{aug}

\received{\smonth{9} \syear{2010}}
\revised{\smonth{4} \syear{2011}}

%
\begin{abstract}
It is known that in many cases distributions of exponential integrals
of L\'evy processes are
infinitely divisible and in some cases they are also selfdecomposable.
In this paper, we give some sufficient conditions under which
distributions of
exponential integrals are not only selfdecomposable but furthermore are
generalized gamma convolution.
We also study exponential integrals of more general independent
increment processes.
Several examples are given for illustration.
\end{abstract}

%

\end{frontmatter}

\section{Introduction}\label{sec1}
Let $(\xi, \eta)=\{(\xi_t,\eta_t), {t\geq0}\}$ be a bivariate c\`
adl\`ag independent
increment process.
In most cases, $(\xi, \eta)$ is assumed as a bivariate L\'evy
process, but we also treat
more general cases where $\xi$ or $\eta$ is a compound sum process,
which is not necessarily a
L\'evy process but is another typical independent increment process.
Our concern in this paper is to examine distributional properties of
the exponential integral
%
\begin{equation}
\label{integral}
V:=\int_{(0,\infty)} \mathrm{e}^{-\xi_{t-}}\,\mathrm{d}\eta_t,
\end{equation}
provided that this integral converges almost surely.
More precisely, we are interested in when $\cL(V)$, the law of $V$,
is selfdecomposable and moreover is a generalized gamma convolution.

We say that a probability distribution $\mu$ on $\R$ (resp. an $\R
$-valued random variable $X$)
is selfdecomposable, if for any $b >1$ there exists a probability
distribution $\mu_b$
(resp. a random variable $Y_b$ independent of $X$) such that
\[
\mu= D_{b^{-1}}(\mu) \ast\mu_b \qquad  (\mbox{resp. }X \overset d =
b^{-1}X+Y_b),
\]
where $D_{a}(\mu)$ means the distribution induced
by $D_{a}(\mu)(aB) := \mu(B)$ for $B\in\mathfrak{B}(\R)$,
$\ast$ is the convolution operator and $\overset d =$ denotes equality
in distribution.
Every selfdecomposable distribution is infinitely divisible.
Some well-known distributional properties of nontrivial selfdecomposable
distributions are absolute continuity and unimodality
(see Sato \cite{Sa99}, pages 181 and 404).

First, we review existing results on $\cL(V)$.
Bertoin \textit{et al.} \cite{BLM06} (in the case when $\eta= \{ \eta_t\}$ is
a one-dimensional L\'evy process) and
Kondo \textit{et al.} \cite{KondoMaejimaSato} (in the case when $\eta$ is a
multi-dimensional
L\'evy process) showed
that if $\xi=\{\xi_t\}$ is a spectrally negative L\'evy process satisfying
$\lim_{t \to\infty} \xi_t = +\infty$ a.s. and if the integral
\eqref{integral} converges a.s., or equivalently, if
$\int_{\R}\log^+ |y| \nu_{\eta}(\mathrm{d}y)<\infty$
for the L\'evy measure $\nu_{\eta}$ of $\eta_1$, then $\cL(V)$ is
selfdecomposable.

On the other hand, there is an example of noninfinitely divisible
$\cL(V)$, which is due to Samorodnitsky (see Kl\"uppelberg \textit{et al.}
\cite{KLM06}).
In fact if $(\xi_t,\eta_t) = (S_t + a t, t)$, where $\{S_t\}$ is a
subordinator and $a>0$ some constant, then the support of $\cL(V)$ is bounded
so that $\cL(V)$ is not infinitely divisible.

Recently, Lindner and Sato \cite{LS09} considered the exponential
integral
\[
\int_{(0,\infty)} \exp ( -(\log c)N_{t-} )\,\mathrm{d}Y_t =\int
_{(0,\infty)}c^{-N_{t-}}\,\mathrm{d}Y_t, \qquad  c>0,
\]
where
$\{(N_t, Y_t)\}$ is a bivariate compound Poisson process whose L\'evy
measure is
concentrated on $(1,0)$, $(0,1)$ and $(1,1)$, and showed a necessary
and sufficient
condition for the infinite divisibility of $\cL(V)$.
They also pointed out that $\cL(V)$ is always $c^{-1}$-decomposable, namely
there exists a probability distribution $\rho$ such that
$\mu= D_{c^{-1}}(\mu) \ast\rho$.
Note that a $c^{-1}$-decomposable distribution is not necessarily
infinitely divisible, unless
$\rho$ is infinitely divisible.
In their second paper (Lindner and Sato \cite{LSpre}), they also gave
a condition under which
$\cL(V)$, generated by a bivariate compound Poisson
process $\{ (N_t, Y_t)\}$ whose L\'evy measure is concentrated on $(1,0)$,
$(0,1)$ and $(1,c^{-1})$, is infinitely divisible.

For other distributional properties of exponential integrals, like the
tail behavior, see, e.g., Maulik and Zwart \cite{MaZw06}, Rivero \cite
{Ripre} and
Behme \cite{behme}.

In this paper, we focus on
\lq\lq Generalized Gamma Convolutions'' (GGCs, for short)
to get more explicit distributional informations of $V$ than
selfdecomposability.

Throughout this paper, we say that for $r>0$ and $\lambda>0$ a random variable
$\gamma_{r,\ld}$ has a $\operatorname{gamma}(r,\ld)$ distribution if its
probability density function
$f$ on $(0,\infty)$ is
\[
f(x) = \frac{\ld^r}{\Gm(r)} x^{r-1}\mathrm{e}^{-\ld x}.
\]
A $\operatorname{gamma}(1,\ld)$ distribution is an exponential distribution with
parameter $\ld>0$.
When we do not have to emphasize the parameters $(r,\ld)$, we just
write $\gamma$
for a gamma random variable.

The class of GGCs is defined to be the smallest class of distributions on
the positive half line that contains all gamma distributions and is
closed under convolution
and weak convergence.
By including gamma distributions on the negative real axis, we obtain
the class
of distributions on $\R$ which will be called \lq\lq Extended Generalized
Gamma Convolutions'' (EGGCs, for short). We refer to Bondesson \cite
{Bondesson} and
Steutel and van~Harn \cite{Steutel} for many properties of GGCs and
EGGCs with
relations among other subclasses of infinitely divisible distributions.

One well-known concrete example of exponential integrals is the following.
When $(\xi_t, \eta_t )= (B_t + a t, t)$ with a standard
Brownian motion $\{B_t\}$ and a drift $a>0$, the law of
\eqref{integral} equals $\cL ({1}/{(2 \gamma)} )$ which
is GGC
(and thus is selfdecomposable).

When choosing $\xi$ to be deterministic, that is, $(\xi_t, \eta_t )=
(t, \eta_t)$,
the exponential integral \eqref{integral} is defined and is an EGGC if and
only if $\eta$ admits a finite log-moment (needed for the convergence) and
$\cL(\eta_1)$ is included in the Goldie--Steutel--Bondesson Class, a superclass
of EGGC as defined for example, in Barndorff-Nielsen \textit{et al.} \cite
{BMS06}. This fact follows directly from
\cite{BMS06}, equation (2.28). In the same paper, the authors
characterized the
class of GGCs by using stochastic integrals with respect to L\'evy processes
as follows. Let
$e(x) = \int_x^{\infty}{u^{-1}}\mathrm{e}^{-u}\,\mathrm{d}u, x>0,$ and let $e^*(\cdot)$
be its
inverse function. Then
$\cL (\int_{(0,\infty)}e^*(s)\,\mathrm{d}\eta_s )$ is GGC
if $\eta$ is a L\'evy process and $\cL(\eta_1)$ has a finite log-moment.

In this paper, via concrete examples, we investigate distributional
properties of
exponential integrals connected with GGCs.

The paper is organized as follows.
In Section~\ref{sec2}, we give some preliminaries. In Section~\ref{sec3}, we consider
exponential integrals for two independent L\'evy processes
$\xi$ and $\eta$
such that $\xi$ or $\eta$ is a compound Poisson process, and
construct concrete examples related to our question.
In the special case that both $\xi$ and $\eta$ are compound Poisson processes,
we also allow dependence between the two components of $(\xi, \eta)$.
In Section~\ref{sec4}, we consider exponential integrals for independent
increment processes such that
$\xi$ and $\eta$ are independent and one is a compound sum process
(which is
not necessarily a L\'evy process) while the other is a L\'evy process.


\section{Preliminaries}\label{sec2}

The class of all infinitely divisible distributions on $\R$
(resp. $\R_+$) is denoted by $I(\R)$ (resp. $I(\R_+)$).
We denote the class of selfdecomposable distributions on $\R$
(resp. $\R_+$) by $L(\R)$ (resp. $L(\R_+)$).
The class of EGGCs on $\R$
(resp. GGCs on $\R_+$) is denoted by $T(\R)$ (resp. $T(\R_+)$).
The moment generating function of a random variable $X$ and of a
distribution $\mu$ are written as $\mathbb{L}_X$ and $\mathbb
{L}_{\mu}$, respectively.
If $X$ is positive and $\mu$ has support in $\R_+$, $\mathbb{L}_X$
and $\mathbb{L}_{\mu}$
coincide with the Laplace transforms.

We are especially interested in distributions on $\R_+$.
The class $T(\R_+)$ is characterized by the Laplace transform as follows:
A probability distribution $\mu$ is GGC if and only if there exist
$a \geq0$ and a measure $U$ satisfying
\[
\int_{(0,1)}|\log x^{-1}|U(\mathrm{d}x)<\infty  \quad \mbox{and} \quad   \int
_{(1,\infty)}x^{-1}U(\mathrm{d}x)<\infty,
\]
such that
the Laplace transform $\mathbb{L}_{\mu}(z)$ can be uniquely
represented as
%
\begin{equation}\label{measureU}
\mathbb{L}_{\mu}(u) = \int_{[0,\infty)} \mathrm{e}^{-ux}\mu(\mathrm{d}x)=
\exp \biggl\{ -au + \int_{(0,\infty)} \log \biggl(\frac
{x}{x+u} \biggr) U(\mathrm{d}x) \biggr\}.
\end{equation}

Another class of distributions which we are interested in is the class of
distributions on $\R_+$ whose densities are hyperbolically completely monotone
(HCM, for short).
Here we say that a function $f(x)$ on $(0,\infty)$ with values in $\R_+$
is HCM if for every $u>0$, the mapping $f(uv)\cdot f(u/v)$ is completely
monotone with respect to the variable $w=v+v^{-1}$, $v>0$.
Examples of HCM functions are $x^{\beta}  (\beta\in\R)$,
$\mathrm{e}^{-cx}$  $(c>0)$
and $(1+cx)^{-\al}$  $(c>0,  \al>0)$.
The class of all distributions on $\R_+$ whose probability densities
are HCM is
denoted by $H(\R_+)$.
Note that $H(\R_+) \subset T(\R_+) \subset L(\R_+) \subset I(\R_+)$.

For illustration, we give some examples.
Log-normal distributions are in $H(\RR_+)$ \cite{Bondesson}, Example
5.2.1. So these are also GGCs. Positive strictly stable
distributions with Laplace transform $\LL(u)=\exp(-u^\alpha)$ for
$\alpha\in\{1/2,1/3,\ldots\}$ are in $H(\RR_+)$ \cite{Bondesson},
Example 5.6.2, while they are GGCs for all $\alpha\in(0,1]$
\cite{Steutel}, Proposition 5.7. Let $Y=\exp(\gamma_{r,\ld})-1$. If
$r \geq1$,
then $\mathcal{L}(Y)$ is in $H(\RR_+)$, but if $r < 1$, $\mathcal
{L}(Y)$ is
not in $H(\RR_+)$ \cite{Bondesson}, page 88. But $\mathcal{L}(Y)$ or
equivalently $\cL(\exp(\gamma_{r,\ld}))$ is always in $T(\RR_+)$,
independent of
the value of $r$ \cite{Bondesson}, Theorem~6.2.3. Remark that by treating
$H(\RR_+)$ we cannot replace $\exp(\gamma_{r,\ld})-1$ by $\exp
(\gamma_{r,\ld})$. Namely, set $r=1$ and observe that the probability density
function
$\lambda(x+1)^{-\lambda-1}1_{[0,\infty)}(x)$ is HCM, but the probability
density function $\lambda x^{-\lambda-1}1_{[1,\infty)}(x)$ is not HCM.
It follows from this that $\mathcal{L}(\exp(\gamma_{1,\ld})-1)$ is in
$H(\RR_+)$ but $\mathcal{L}(\exp(\gamma_{1,\ld}))$ is not in
$H(\RR_+)$.


In addition, we also investigate the modified HCM class denoted by
$\wt H(\R)$, which gives some interesting examples of $\cL(V)$ on $\R$.
The class $\wt H(\R)$ is characterized to be the class of
distributions of random
variables $\sqrt{X}Z$, where $\cL(X)\in H(\R_+)$ and $Z$ is a
standard normal
random variable independent of $X$ (see Bondesson \cite{Bondesson}, page 115).
By the definition, any distribution in $\wt H(\R)$ is a type $G$
distribution, which is the
distribution of the variance mixture of a standard normal random variable.
Note that    $\wt H(\R) \subset T(\R)$.
As will be seen in Proposition~\ref{facts}, there are nice relations between
$\wt H(\R)$ and $T(\R)$ in common with those of $H(\R_+)$ and $T(\R_+)$.

Here we state some known facts that we will use later.

\begin{prop}[(Bondesson \cite{Bondesson} and Steutel and van~Harn
\cite{Steutel})]
\label{facts}
\begin{longlist}[(5)]
\item[(1)] A continuous function $\LL(u),  u>0,$
with $\LL(0+)=1$ is HCM if and only if it is the Laplace transform of
a GGC.

\item[(2)] If $\cL(X)\in H(\R_+)$, $\cL(Y)\in T(\R_+)$ and $X$ and $Y$ are
independent, then $\cL(XY)\in T(\R_+)$.

\item[(3)] Suppose that $\cL(X) \in H(\R_+)$, $\cL(Y) \in T(\R)$ and that
$X$ and $Y$ are independent.
If $\cL(Y)$ is symmetric, then $\cL(\sqrt{X}Y)\in T(\R)$.

\item[(4)] Suppose that $\cL(X) \in\wt H(\R)$ and $\cL(Y) \in T(\R)$
and that
$X$ and $Y$ are independent.
If $\cL(Y)$ is symmetric, then $\cL(X Y) \in T(\R)$.

\item[(5)] If $\cL(X)\in\wt H(\R)$, then $\cL(|X|^q)\in H(\R_+)$ for
all $|q|\ge2,  q\in\R$.
Furthermore, $\cL(|X|^q \*\operatorname{sign} (X)) \in\wt H(\R)$
for all $q\in\mathbb{N}, q\neq2$, but not always for $q=2$.
\end{longlist}
\end{prop}

\begin{rem}
Notice that the distribution
of a sum of independent random variables with distributions in $H(\R_+)$
does not necessarily belong to $H(\R_+)$. See Bondesson \cite
{Bondesson}, page 101.
\end{rem}

Some distributional properties of GGCs are stated in the
following proposition \cite{Bondesson}, Theorems~4.1.1. and 4.1.3.

\begin{prop}
\begin{longlist}[(2)]
\item[(1)] The probability density function of a GGC without Gaussian part
satisfying $0<\int_{(0,\infty)} U(\mathrm{d}u)=\beta<\infty$ with the
measure $U$ as in \eqref{measureU}
admits the representation $x^{\beta-1}h(x)$,
where $h(x)$ is some completely monotone function.

\item[(2)] Let $f$ be the probability density of a GGC distribution without
Gaussian part
satisfying $1< \int_{(0,\infty)} U(\mathrm{d}u)=\beta\leq\infty$.
Let $k$ be a nonnegative integer such that $k < \beta-1$.
Then the density $f$ is continuously differentiable any times on
$(0,\infty)$,
and at $0$ at least
$k$ times differentiable with $f^{(j)}(0)=0$ for $j \leq k$.
\end{longlist}
\end{prop}

Examples of GGCs and the explicit calculation of their L\'evy measure
are found in Bondesson \cite{Bondesson} and James \textit{et al.} \cite{JRY08}.

Necessary and sufficient conditions for the convergence of \eqref
{integral} for
bivariate L\'evy processes were given by Erickson and Maller \cite
{EricksonMaller}.
More precisely, in their Theorem 2.1, they showed that $V$ converges
a.s. if and only if for some
$\epsilon>0$ such that $A_\xi(x)>0$ for all $x>\epsilon$ it holds
%
\begin{eqnarray} \label{cond-of-conv}
 &\displaystyle \lim_{t\to\infty}\xi_t=\infty   \qquad \mbox{a.s.} \quad  \mbox{and}&
\nonumber
\\[-8pt]
\\[-8pt]
&\displaystyle \int_{(\mathrm{e}^{\epsilon},\infty)} \biggl(\frac{\log y}{A_\xi(\log
y)} \biggr)
|\nu_\eta((\mathrm{d}y,\infty))+\nu_\eta((-\infty,-\mathrm{d}y))| <\infty.&
\nonumber
\end{eqnarray}
Here $A_\xi(x)=a_\xi+ \nu_\xi((1,\infty))
+\int_{(1,x]}\nu_\xi((y,\infty))\,\mathrm{d}y$ while $(\Sigma_X, \nu_X,
a_X)$ denotes
the L\'evy--Khintchine triplet of a L\'{e}vy process $X$.


\section{Exponential integrals for compound Poisson processes}\label{sec3}

In this section, we study exponential integrals of the form
\eqref{integral}, where either $\xi$ or $\eta$ is a compound
Poisson process and the other is an arbitrary L\'evy process.
First, we assume the two processes to be independent, later we
also investigate the case that $(\xi,\eta)$ is a bivariate
compound Poisson process.

\subsection{Independent component case}\label{sec3.1}

We start with a general lemma which gives a sufficient condition for
distributions of perpetuities to be GGCs.

\begin{lem} \label{perpetuityinthorin}
Suppose $A$ and $B$ are two independent random
variables such that $\cL(A)\in H(\R_+)$ and $\cL(B)\in T(\R_+)$.
Let $(A_j,B_j), j=0,1,2,\ldots,$ be i.i.d. copies of $(A,B)$.
Then, given its a.s. convergence, the distribution of the perpetuity
$Z:=\sum_{k=0}^\infty ( \prod_{i=0}^{k-1} A_i ) B_k$
belongs to
$T(\R_+)$. Furthermore, if $\cL(A)\in\wt H(\R)$, $\cL(B)\in T(\R
)$ and
$\cL(B)$ is symmetric, then $\cL(Z)\in T(\R)$.
\end{lem}

\begin{pf}
If we put
\[
Z_n:= \sum_{k=0}^n  \Biggl(\prod_{i=0}^{k-1} A_i \Biggr) B_k,
\]
then we can rewrite
\[
Z_n=B_0+ A_0\bigl(B_1 + A_1\bigl(B_2 +\cdots+ A_{n-2}(B_{n-1}+A_{n-1}B_n) \cdots\bigr)\bigr).
\]
Since $A_{n-1}$ and $B_n$ are independent, $\cL(A_{n-1})\in H(\R_+)$ and
$\cL(B_n) \in T(\R_+)$,
we get $\cL(A_{n-1}B_n)\in T(\R_+)$ by Proposition~\ref{facts}.
Further $B_{n-1}$ and $A_{n-1}B_n$ are independent and both
distributions belong to $T(\R_+)$
and hence $\cL(B_{n-1}+A_{n-1}B_n)\in T(\R_+)$.
By induction, we can conclude that $\cL(Z_n)\in T(\R_+)$.
Since the class $T(\R_+)$ is closed under weak convergence,
the first part follows immediately.
A similar argument with $(4)$ in Proposition~\ref{facts} gives the
second part.
\end{pf}

\subsubsection*{Case 1: The process $\xi$ is a compound Poisson process}

\begin{prop}\label{CPPinThorin}
Suppose that the processes $\xi$ and $\eta$ are
independent L\'evy processes and further that $\xi_t=\sum
_{i=1}^{N_t}X_i$ is a
compound Poisson process
with i.i.d. jump heights $X_i, i=1,2,\ldots,$ such that
$0<E[X_1]<\infty$,
$\cL(\mathrm{e}^{-X_1})\in H(\R_+)$ and $\eta$ has finite log-moment $E\log
^+|\eta_1|$
and it holds $\cL(\eta_\tau) \in T(\R_+)$ for an
exponential random variable $\tau$ independent of $\eta$ having the same
distribution as the waiting times of $N$.
Then the integral \eqref{integral} converges a.s.  and
\[\cL\biggl (\int_{(0,\infty)} \mathrm{e}^{-\xi_{t-}}\,\mathrm{d}\eta
_t \biggr) \in
T(\R_+).
\]
Furthermore, if $\cL(\mathrm{e}^{-X_1})\in\wt H(\R)$, $\cL(\eta_\tau)\in
T(\R)$ and
$\cL(\eta_\tau)$ is symmetric, then $\cL(V)\in T(\R)$.
\end{prop}

\begin{pf}
Convergence of the integral follows from \eqref{cond-of-conv}.

Set $T_0=0$ and let $T_j, j=1,2,\ldots,$ be the time of the $j$th jump of
$\{N_t,t\geq0\}$.
Then we can write
\begin{eqnarray*}
\int_{(0,\infty)} \mathrm{e}^{-\xi_{t-}}\,\mathrm{d}\eta_t &=& \sum_{j=0}^\infty\int
_{(T_j,T_{j+1}]}
\mathrm{e}^{-\sum_{i=1}^j X_i}\,\mathrm{d}\eta_t
= \sum_{j=0}^\infty \mathrm{e}^{-\sum_{i=1}^j X_i}
\int_{(T_j,T_{j+1}]}\,\mathrm{d}\eta_t \\
&\hspace*{2.3pt}=:& \sum_{j=0}^\infty \Biggl(\prod_{i=1}^j A_i  \Biggr) B_j,
\end{eqnarray*}
where $\sum_{i=1}^0=0$, $\prod^0_{i=1}=1$, $A_i=\mathrm{e}^{-X_i}$ and $B_j=\int_{(T_j,T_{j+1}]}\,\mathrm{d}\eta_t \overset
{d}=\eta_{_{ T_{j+1}-T_j}}$.
Now Lemma~\ref{perpetuityinthorin} yields the conclusion.\vspace*{-2pt}
\end{pf}

In the following, we first give some examples for possible choices of
$\xi$
fulfilling the assumptions of Proposition~\ref{CPPinThorin} and then continue
with examples for $\eta$. Hence, any combination of them yields an exponential
integral which is a GGC.\vspace*{-2pt}

\begin{ex}[(The case when $X_1$ is a normal random variable with
positive mean)]
We see that $\cL(\mathrm{e}^{-X_1})$ is log-normal and hence is in $H(\R_+)$.\vspace*{-2pt}
\end{ex}

\begin{ex}[(The case when $X_1$ is the logarithm of the power of a gamma random
variable)]
Let $Y = \gamma_{r,\ld}$ and $X_1=c \log Y$ for $c\in\mathbb{R}$.
Recall that $\cL(\log Y)\in T(\R)$ and so
$\cL(c X_1)=\cL(\log Y^c)\in T(\R)$ for $c\in\R$.
Note that $E[X_1]=c \psi(\lambda)$,
where $\psi(x)$ the derivative of $\log\Gamma(x)$.
If we take $c\in\R$ such that $c \psi(\lambda)>0$, we conclude that
$\cL(\mathrm{e}^{-X_1})=\cL(\gamma_{r,\ld}^{-c})\in H(\R_+)$.\vspace*{-2pt}
\end{ex}

\begin{ex}[(The case when $X_1$ is logarithm of some positive strictly
stable random variable)]
Let $X_1=\log Y$ be a random variable,
where $Y$ is a positive stable random variable with parameter $0<\alpha<1$.
Then $X_1$ is in the class of EGGCs when
$\alpha=1/n, n=2,3,\ldots$ (see Bondesson \cite{Bondesson}, Example
7.2.5) and
\[
E[\mathrm{e}^{uX_1}]=E[Y^u]=\frac{\Gamma(1-u/\alpha)}{\Gamma(1-u)}.
\]
It follows that
$E[X_1]=-\frac{1}{\alpha}\psi(1)+\psi(1) = (1-1/\alpha)\psi(1)>0$
and $\cL(\mathrm{e}^{-X_1})=\cL(Y^{-1})\in H(\RR_+)$ by \cite{Bondesson},
Example 5.6.2.\vspace*{-2pt}
\end{ex}

\begin{ex}[(The case when $X_1$ is the logarithm of the ratio of two
exponential random variables)]
Let $X_1=\log(Y_1/Y_2)$, where $Y_j, j=1,2,$ are independent exponential
random variables with parameters $\ld_j >0,  j=1,2$.
The density function of $X_1$ is given by \cite{Bondesson}, Example
7.2.4,
\[
f(x)=\frac{1}{B(\lambda_1,\lambda_2)}\frac{\mathrm{e}^{-\lambda_1
x}}{(1+\mathrm{e}^{-x})^{\lambda_1+\lambda_2}},
 \qquad   x\in\mathbb{R},
\]
where $B(\cdot,\cdot)$ denotes the Beta-function.
Now if $E[X_1]>0$ we can set $X_1$ to be a jump distribution of $\xi$.
It is easy to see that $\cL(\mathrm{e}^{-X_1}) =
\cL(Y_2/Y_1)\in H(\mathbb{R}_+)$ since $\cL(Y_j)\in H(\mathbb{R}_+)$.\vspace*{-2pt}
\end{ex}

\begin{ex}[(The case when $\eta$ is nonrandom)]\label{ex1}
When $\eta_t=t$, it holds that $\cL(\eta_\tau)=\cL(\tau)\in T(\R_+)$.\vspace*{-2pt}
\end{ex}

\begin{ex}[(The case when $\eta$ is a stable subordinator)]\label{ex2}
The Laplace transform of
$B:=\eta_\tau$ with $\tau=\gamma_{1,\ld}$ is given by (see, e.g.,
Steutel and van~Harn \cite{Steutel}, page 10)
\[
\mathbb{L}_B(u)=\frac{\lambda}{\lambda- \log\mathbb{L}_\eta(u)}.\vadjust{\goodbreak}
\]
Now consider $\eta$ to be a stable subordinator without drift. Then the
Laplace transform of $\eta_1$ is given by $\mathbb{L}_{\eta_1} (u) =
\exp
\{-u^{\alpha}\}$. Therefore, we have
\[
\mathbb{L}_B(u)=\frac{\lambda}{\lambda+ u^{\alpha} }.
\]
This function is HCM, since $\frac{\lambda}{\lambda+u}$ is HCM by the
definition and due to the fact that the composition of an HCM function and
$x^\alpha$, $|\alpha|\leq1$, is also HCM (see \cite{Bondesson},
page 68).
Thus, the Laplace transform of $B$ is HCM by Proposition~\ref{facts}
and we
conclude that $\cL(\eta_\tau)$ is GGC.

Remark that if $\eta$ admits an additional drift term, the distribution
$\cL(B)$ is not GGC. This result was pointed out by Kozubowski \cite{Ko05}.
\end{ex}

\begin{ex}[(The case when $\eta$ is an inverse Gaussian L\'evy process)]
We suppose $\eta$ to be an inverse Gaussian subordinator with
parameters $\beta>0$ and $\delta>0$. The
Laplace transform of $\eta_t$ is
\[
\mathbb{L}_{\eta_t} (u)= \exp \bigl( -\delta t\bigl(
\sqrt{\beta^2 +2u}-\beta\bigr) \bigr).
\]
Now by choosing the parameters satisfying $\lambda\geq\delta\beta$, we
have, for $B=\eta_{\tau}$,
\[
\mathbb{L}_B(u)=\frac{\lambda}{\lambda-\delta\beta+ \delta
\sqrt{\beta^2+2u}}.
\]
This function is HCM by argumentation as in Example~\ref{ex2} with
$\alpha=1/2$
and using Property (xi) in \cite{Bondesson}, page 68.
\end{ex}

\begin{rem}
Although the L\'evy measure of $\eta_\tau$ is known explicitly, it is
an open
question whether the parameter of the exponentially distributed random
variable $\tau$ has an influence on the GGC-property of $\eta_\tau$, or
not. Examples~\ref{ex1} and~\ref{ex2} lead to the conjecture that
there is no
influence. So far, no counterexamples to this conjecture are known to the
authors.
\end{rem}

\subsubsection*{Case 2: The process $\eta$ is a compound Poisson process}
In the following, we assume the integrator $\eta$ to be a compound Poisson
process, while $\xi$ is an arbitrary L\'evy process, independent of
$\eta$.
We can argue similarly as above to obtain the following result.

\begin{prop}\label{CPPinThorin2}
Let $\xi$ and $\eta$ be independent and assume
$\eta_t=\sum_{i=1}^{N_t}Y_i$ to be a compound Poisson process with
i.i.d. jump
heights $Y_i, i=1,2,\ldots.$
Suppose that $E[\xi_1]>0$, $E\log^+|\eta_1|<\infty$, $\cL(Y_1)\in
T(\R_+)$
and $\cL(\mathrm{e}^{-\xi_\tau})\in H(\R_+)$ for an exponentially
distributed random
variable $\tau$ independent of $\xi$ having the same distribution as the
waiting times of~$N$.
Then the integral \eqref{integral} converges a.s. and it holds that
\[\cL \biggl(\int_{(0,\infty)} \mathrm{e}^{-\xi_{t-}}\,\mathrm{d}\eta
_t \biggr) \in
T(\R_+).
\]
Furthermore, if $\cL(\mathrm{e}^{-\xi_\tau})\in\wt H(\R)$, $\cL(Y_1)\in
T(\R)$ and
$\cL(Y_1)$ is symmetric, then $\cL(V)\in T(\R)$.
\end{prop}

\begin{pf}
Convergence of the integral is guaranteed by \eqref{cond-of-conv}.
Now set $T_0=0$ and let $T_j, j=1,2,\ldots,$ be the jump times of $\{
N_t,t\geq0\}$.
Then we have
\begin{eqnarray*}
\int_{(0,\infty)} \mathrm{e}^{-\xi_{t-}}\,\mathrm{d}\eta_t
&=& \sum_{j=1}^\infty \mathrm{e}^{-\xi_{T_j}} Y_j
= \sum_{j=1}^\infty \mathrm{e}^{-(\xi_{T_j}-\xi_{T_{j-1}})}\cdots
\mathrm{e}^{-(\xi_{T_1}-\xi_{T_0})} Y_j \\
&=& \sum_{j=1}^\infty \Biggl(\prod_{i=1}^j \mathrm{e}^{-(\xi_{T_i}-\xi
_{T_{i-1}})} \Biggr)
Y_j
=: \sum_{j=1}^\infty \Biggl(\prod_{i=1}^j A_i  \Biggr) B_j,
\end{eqnarray*}
where $A_i=\mathrm{e}^{-(\xi_{T_i}-\xi_{T_{i-1}})}\overset{d}= \mathrm{e}^{-\xi
_{T_i-T_{i-1}}}$
and $B_j=Y_j$.
Remark that the proof of Lemma~\ref{perpetuityinthorin}
remains valid even if the summation starts from $j=1$.
Hence, the assertion follows from Lemma~\ref{perpetuityinthorin}.
\end{pf}


\subsection{Dependent component case}\label{sec3.2}
In this subsection, we generalize a model of Lindner and Sato \cite
{LS09} and study to which
class $\cL(V)$ belongs.

Let $0<p<1$.
Suppose that $(\xi, \eta)$ is a bivariate compound Poisson process
with parameter $\lambda>0$ and normalized L\'evy measure
\[
\nu(\mathrm{d}x,\mathrm{d}y) = p \delta_0(\mathrm{d}x)\rho_0(\mathrm{d}y) + (1-p) \delta_1(\mathrm{d}x) \rho_1(\mathrm{d}y),
\]
where $\rho_0$ and $\rho_1$ are probability measures on $(0,\infty)$ and
$[0,\infty)$, respectively, such that
\[
\int_{(1,\infty)}\log y\,\mathrm{d}\rho_0(y)<\infty  \quad \mbox{and} \quad  \int
_{(1,\infty)}\log y\,\mathrm{d}\rho_1(y)<\infty.
\]
For the bivariate compound Poisson process $(\xi, \eta)$, we have the
following representation
(see Sato \cite{Sa99}, page 18):
\[
 ( \xi_t , \eta_t  )
= \sum_{k=0}^{N_t} {S}_{k}
=  \Biggl(\sum_{k=0}^{N_t} S^{(1)}_{k},
\sum_{k=0}^{N_t} S^{(2)}_{k}  \Biggr),
\]
where $S^{(1)}_{0}=S^{(2)}_{0}=0$ and $ \{{S}_{k} =
 ( S^{(1)}_k , S^{(2)}_k  )
 \}_{k=1}^{\infty}$ is a sequence of two-dimensional i.i.d.
random variables.
It implies that the projections of the compound Poisson process on $\R^2$
are also compound Poisson processes.
Precisely in the given model, since
$\mathbb{P}(S^{(1)}_1 = 0)=p$ and $\mathbb{P}(S^{(1)}_1 = 1)=1-p$,
the marginal
process $\xi$ is a Poisson process with parameter $(1-p)\lambda>0$.
Note that $S^{(1)}_{k}$ and $S^{(2)}_{k}$ may be dependent for any $k
\in\N$.
In this case, $\rho_i (B)$ is equal to $\mathbb{P}(S^{(2)}_k \in B |
S^{(1)}_k =i)$
for $i=0,1$ and $B \in\mathcal{B}(\R)$.

\begin{ex}\label{LSmodel}
In Lindner and Sato \cite{LS09}, the authors considered the bivariate
compound Poisson process
with parameter $u+v+w$, $u,v,w\geq0$ and normalized L\'evy measure
\[
\nu(\mathrm{d}x,\mathrm{d}y) = \frac{v}{u+v+w}\delta_0(\mathrm{d}x)\delta_1(\mathrm{d}y)
+ \frac{u+w}{u+v+w}\delta_1(\mathrm{d}x) \biggl(
\frac{u}{u+w}\delta_0(\mathrm{d}y) + \frac{w}{u+w}\delta_1(\mathrm{d}y) \biggr).
\]
So in their setting $p= \frac{v}{u+v+w}$, $\rho_0 = \delta_1$
and $\rho_1= \frac{u}{u+w}\delta_0 + \frac{w}{u+w}\delta_1$.
\end{ex}

In the following theorem, we give a sufficient condition for $\cL(V)$,
given by
\eqref{integral} with $(\xi,\eta)$ as described above, to be GGC.

\begin{thm}\label{main3.2}
If the function
$\frac{(1-p)\mathbb{L}_{\rho_1} (u)}{1-p\mathbb{L}_{\rho_0} (u)}$
is HCM,
then $\cL(V)$ is GGC.
\end{thm}

\begin{pf}
Define $T_\xi$ and $M$ to be the first jump time of the Poisson
process $\xi$ and the
number of the jumps of the bivariate compound Poisson process
$(\xi,\eta)$ before $T_\xi$, respectively.
Due to the strong Markov property of the L\'evy process $(\xi,\eta)$,
we have
\begin{eqnarray}
\int_{(0,\infty)} \exp(-\xi_{s-})\,\mathrm{d}\eta_s
&=& \int_{(0,T_{\xi}]} \exp(-\xi_{s-})\,\mathrm{d}\eta_s + \int_{(T_{\xi
},\infty]}
\exp(-\xi_{s-})\,\mathrm{d}\eta_s\nonumber\\
&=& \eta_{T_{\xi}} + \int_{(0,\infty)} \exp(-\xi_{T_{\xi}+s-})\,\mathrm{d}\eta_{s+T_{\xi}}\nonumber\\
&=& \eta_{T_{\xi}} + \exp(-\xi_{T_{\xi}})\int_{(0,\infty)}
\exp\bigl(-(\xi_{T_{\xi}+s-} -\xi_{T_{\xi}})\bigr) \,\mathrm{d}\bigl(\bigl(\eta_{(s+T_{\xi})} -
\eta_{_{T_{\xi}}}\bigr)
+\eta_{_{T_{\xi}}}\bigr)\nonumber\\
&\overset{d}{=}& \eta_{_{T_{\xi}}}
+ \mathrm{e}^{-1}\int_{(0,\infty)} \exp(-\widetilde{\xi}_{s-})\,\mathrm{d}\widetilde{\eta}_{s},\nonumber
\end{eqnarray}
where the process $(\widetilde{\xi},\widetilde{\eta})$ is
independent of
$\{(\xi_t,\eta_t),t\leq T_\xi\}$ and has the same law as $(\xi,\eta)$.
Therefore, we have
%
\begin{equation}\label{ssd}
\mathbb{L}_{\mu}(u) =\mathbb{L}_{\mu}(\mathrm{e}^{-1}u) \mathbb{L}_{\rho}(u),
\end{equation}
with $\mu=\cL(V)$ and $\rho$ denoting the distribution of
$\eta_{_{T_{\xi}}}$. Thus, $\mu$ is $\mathrm{e}^{-1}$-decomposable and it
follows that
\[
\mathbb{L}_{\mu}(u) = \prod_{n=0}^{\infty} \mathbb{L}_{\rho}(\mathrm{e}^{-n}u).
\]
In the given setting, we have
\begin{eqnarray*}
\eta_{_{T_\xi}} = \Biggl (\sum_{k=0}^{M} S^{(2)}_k  \Biggr) + S^{(2)}_{M+1},
\end{eqnarray*}
where $M$ is geometrically distributed with parameter $p$, namely,
\[
\mathbb{P}(M=k) = (1-p)p^{k}   \qquad \mbox{for any }k \in\mathbb{N}_0.
\]
Hence, we obtain
\begin{eqnarray*}
\mathbb{L}_{\rho} (u) &=& \mathbb{E} [\exp(-u \eta_{_{T_\xi
}}) ] \\
&=& \mathbb{E} \Biggl[\mathbb{E} \Biggl[\exp \Biggl(-u \Biggl(\sum
_{k=0}^{M} S^{(2)}_k
+ S^{(2)}_{M+1} \Biggr) \Biggr) \Big|M \Biggr] \Biggr]= \mathbb{E} [(\mathbb{L}_{\rho_0} (u))^{M}\mathbb{L}_{\rho
_1} (u) ]\\
&=& (1-p)\mathbb{L}_{\rho_1} (u) \sum_{k=0}^{\infty}(p\mathbb
{L}_{\rho_0} (u))^{k}=\frac{(1-p)\mathbb{L}_{\rho_1} (u)}{1-p\mathbb{L}_{\rho_0} (u)}.
\end{eqnarray*}
The class of HCM functions is closed under scale transformation,
multiplication and limit. Therefore, $\mathbb{L}_{\mu}(u)$ is HCM if
$\mathbb{L}_{\rho}(u)$ is HCM and hence $\mu$ is GGC if $\rho$ is
GGC by
Proposition~\ref{facts}(6). As a result, $\mu$ is GGC if
$\frac{(1-p)\mathbb{L}_{\rho_1} (u)} {1-p\mathbb{L}_{\rho_0} (u)}$
is HCM.
\end{pf}

A distribution
with the Laplace transform $\frac{(1-p)}{1-p\mathbb{L}_{\rho_0} (u)}$
is called a compound geometric distribution.
It is compound Poisson, because every
geometric distribution is compound Poisson with L\'evy measure given by
\[
\nu_p(\{ k\}) = -\frac{1}{\log(1-p)} \frac{1}{k+1}p^{k+1}, \qquad
k=1,2,\ldots
\]
(see page 147 in Steutel and van~Harn \cite{Steutel}). Since HCM
functions are closed under
multiplication, we have the following.

\begin{cor} If $\rho_1$ and the compound geometric distribution of
$\rho_0$
are GGCs, then so is $\cL(V)$.
\end{cor}

In addition, we observe the following.

\begin{cor}
\begin{longlist}[(3)]
\item[(1)] For any $c>1$, the distribution
$\mu_c = \law ( \int_{(0,\infty)} c^{-\xi_{s-}}\,\mathrm{d}\eta_s )$
is $c^{-1}$-selfdecomposable. Thus, in the nondegenerate case it is absolutely
continuous or continuous singular (Wolfe~\cite{Wolfe}) and Theorem
\ref{main3.2} holds true also for $\mu_c$ instead of $\mu$.

\item[(2)] If $\rho_1$ is infinitely divisible, then $\mu_c$ is also infinitely
divisible.

\item[(3)] Let $B(\R_+)$ be
the Goldie--Steutel--Bondesson class, which is
the smallest class that contains all mixtures of exponential distributions
and is closed under convolution and weak convergence.

If
$\frac{(1-p)}{1-p\mathbb{L}_{\rho_0} (u)}$
is the Laplace transform of a distribution in $B(\R_+)$, then $\mu_c$
is in $B(\R_+)$.
Moreover, $\mu_c$ will be a $c^{-1}$-semi-selfdecomposable distribution.
\end{longlist}
\end{cor}

About the definition and basic properties of semi-selfdecomposable
distributions, see~\cite{Sa99}. The proof of (1) is obvious.
For (2), remark that a distribution with Laplace
transform $\frac{(1-p)}{1-p\mathbb{L}_{\rho_0} (u)}$ as compound Poisson
distribution is always infinitely divisible.
Hence, only $\rho_1$ has influence
on that property.
The proof of (3) follows from the characterization of the class
$B(\R_+)$ in Chapter $9$ of \cite{Bondesson} and our proof of Theorem
\ref{main3.2}.

\begin{ex}\label{expex}
Let $\rho_1$ be a GGC, that is, $\mathbb{L}_{\rho_1} (u)$ is HCM.
Then if
$\frac{(1-p)}{1-p\mathbb{L}_{\rho_0} (u)}$ is HCM, $\mu$ is found
to be GGC.
For example, if $\rho_0$ is an
exponential random variable with density $f(x)=b\mathrm{e}^{-bx}, b>0$, then
$\frac{1}{1-p\mathbb{L}_{\rho_0} (u)}$ is HCM.
To see this, for $u>0,
v>0$, write
\begin{eqnarray*}
\frac{1}{1-p\mathbb{L}_{\rho_0} (uv)} \frac{1}{1-p\mathbb{L}_{\rho
_0}(u/v)}
&= & \frac{1+\fracd{u}{b}(v+v^{-1})+\fraca{u^2}{b^2}}{(1-p)^2+\fracd
{u}{b}(v+v^{-1})(1-p)+\fraca{u^2}{b^2}}\\
& =& \frac{1}{1-p} + \frac{p+ (1-\fracc{1}{1-p} )\fraca
{u^2}{b^2}}{(1-p)^2+\fracd{u}{b}(v+v^{-1})(1-p)+\fraca{u^2}{b^2}}.
\end{eqnarray*}
This is nonnegative and completely monotone as a function of $v+v^{-1}$.
\end{ex}

\begin{ex}
In the case of Example~\ref{LSmodel}, the L\'evy measure of $\mu_c$ is
\[
\nu_{\mu_c}=\sum_{n=0}^{\infty}\sum_{m=1}^{\infty} a_m \delta_{c^{-n}m},
\]
where $a_m = \frac{1}{m}(q^m-(-r/p)^m)$.
This L\'evy measure is not absolutely continuous.
Thus, $\mu_c$ is never GGC for any parameters $u$, $v$, $w$ and $c$.
\end{ex}


\section{Exponential integrals for independent increment processes}\label{sec4}

We say that a process $X=\{X_t=\sum_{i=1}^{M_t} X_i, t\ge0\}$ is a compound
sum process, if the
$\{ X_i\}$ are i.i.d. random variables, $\{M_t, t\geq0\}$ is a
renewal process and they are independent.
When $\{M_t\}$ is a Poisson process, $X$ is nothing but a compound
Poisson process and
is a L\'evy process.
Unless $\{M_t\}$ is a Poisson process, $X$ is no L\'evy process.
In this section, we consider the case when either $\xi$ or $\eta$
is a compound sum process and the other is an arbitrary L\'evy process.
Although $(\xi, \eta)$ is not a L\'evy process,
the exponential integral \eqref{integral} can be defined and its
distribution can be infinitely divisible and/or GGC in many cases as we
will show
in the following.

\subsection*{Case 1: The process $\xi$ is a compound sum process}
First, we give a condition for the convergence of the exponential
integral \eqref{integral} when $(\xi, \eta)$ is not a L\'evy process.

\begin{prop} \label{convcondxirenewal}
Suppose that $(\xi_t,\eta_t)_{t\geq0}$ is a stochastic process where
$\xi$
and $\eta$ are independent, $\eta$ is a L\'evy process and $\xi
_t=\sum_{i=1}^{M_t}X_i$ is a compound sum process with i.i.d. jump
heights $X_i, i=1,2,\ldots,$ and i.i.d. waiting times $W_i$. Then
\eqref{integral} converges in probability to a finite random variable
if and only if
%
\begin{equation} \label{eq-convcondxirenewal}
\xi_t\to\infty  \qquad \mbox{a.s.}   \quad \mbox{and} \quad    \int_{(1,\infty)}
 \biggl( \frac{\log q}{A_\xi(\log q)}
 \biggr) P(|\eta_{_{W_1}}|\in \mathrm{d}q)<\infty
\end{equation}
for $A_\xi(x)=\int_{(0,x)} P(X_1>u)\,\mathrm{d}u.$
\end{prop}

\begin{pf}
As argued in the proof of Proposition~\ref{CPPinThorin}, we can
rewrite the
exponential integral as perpetuity\vspace*{-2pt}
\[
\int_{(0,\infty)} \mathrm{e}^{-\xi_{t-}}\,\mathrm{d}\eta_t = \sum_{j=0}^\infty
\Biggl(\prod_{i=1}^j A_i  \Biggr) B_j,
\]
where $A_i=\mathrm{e}^{-X_i}$ and $B_j \overset{d}= \eta_{_{W_j}}$.
By Theorem 2.1 of
\cite{goldiemaller00} the above converges a.s. to a
finite random variable if and only if $\prod_{i=1}^n A_i \to0$ a.s. and
\[
\int_{(1,\infty)}  \biggl( \frac{\log q}{A(\log q)}
 \biggr) P(|B_1|\in \mathrm{d}q)<\infty
\]
for $A(x)=\int_{(0,x)}P(-\log A_1>u)\,\mathrm{d}u$.
Using the given expressions for $A_1$ and $B_1$ in our setting, we
observe that
this is equivalent to \eqref{eq-convcondxirenewal}. It remains to show that
a.s. convergence of the perpetuity implies convergence in probability of
\eqref{integral}. Therefore, remark that
\[
\int_{(0,t]}\mathrm{e}^{-\xi_{s-}}\,\mathrm{d}\eta_s=\int_{(0,T_{M_t}]}\mathrm{e}^{-\xi
_{s-}}\,\mathrm{d}\eta_s +
\mathrm{e}^{-\xi_{T_{M_t}}}(\eta_t-\eta_{_{T_{M_t}}}),
\]
where the first term converges to a finite random variable while the second
converges in probability to $0$ since $\sup_{t\in
[T_{M_t},T_{M_{t+1}})}|\eta_t-\eta_{_{T_{M_t}}}|\overset{d}=\sup_{t\in
[0,W_1)}|\eta_t|$.
\end{pf}

Now we can extend Proposition~\ref{CPPinThorin} in this new setting as follows.

\begin{prop}\label{renewalinThorin}
Suppose that the processes $\xi$ and $\eta$ are
independent and that $\xi_t=\sum_{i=1}^{M_t}X_i$ is a compound sum
process with
i.i.d. jump heights $X_i, i=1,2,\ldots,$ and i.i.d. waiting times
$W_i,  i=1,2,\ldots,$ such that
\eqref{eq-convcondxirenewal} is fulfilled.
Suppose that $\cL(\mathrm{e}^{-X_1})\in H(\R_+)$ and
$\cL(\eta_\tau) \in T(\R_+)$ for $\tau$ being a random variable
with the same distribution
as $W_1$ and independent of $\eta$.
Then
\[
\cL\biggl (\int_{(0,\infty)} \mathrm{e}^{-\xi_{t-}}\,\mathrm{d}\eta_t \biggr) \in T(\R_+).
\]
Furthermore, if $\cL(\mathrm{e}^{-X_1})\in\wt H(\R)$, $\cL(\eta_\tau)\in
T(\R)$ and
$\cL(\eta_\tau)$ is symmetric, then $\cL(V)\in T(\R)$.
\end{prop}

In the following, we give some examples fulfilling the assumptions of
Proposition~\ref{renewalinThorin}.

\begin{ex}[(The case when $\eta$ is nonrandom and $\cL(W_1)$ is GGC)]
For the case $\eta_t=t$, $\cL(\eta_\tau)$ belongs to $T(\R_+)$ if
and only if $\cL(\tau)$ does.
Hence, for all waiting times which are GGCs and for a suitable jump
heights of $\xi$,
we have $\cL(V)\in T(\R_+).$
\end{ex}

\begin{ex}[(The case when $\eta$ is a stable subordinator and $\cL
(W_1)$ is GGC)]
Consider $\eta$ to be a stable subordinator having Laplace transform
$\mathbb{L}_{\eta} (u) =
\exp\{-u^{\alpha}\}$ with $0<\alpha< 1$.
Then the Laplace transform of $B:=\eta_\tau$ is given by
$
\mathbb{L}_B(u)=\LL_{\tau} (u^{\alpha} ).
$
This function is HCM if and only if $\tau$ is GGC, since by
Proposition~\ref{facts},
$\LL_\tau$ is HCM and hence also its composition with $x^\alpha$.
Thus, whenever $\cL(\tau)=\cL(W_1)$ is GGC, $\cL(\eta_\tau)$ is
GGC, too,
fulfilling the assumption of Proposition~\ref{renewalinThorin}.
\end{ex}

\begin{ex}[(The case when $\eta$ is a standard Brownian motion and $\cL(W_1)$
is GGC)]
Given that $\eta$ is a standard Brownian motion, $\eta_1$ has
characteristic function $E\mathrm{e}^{\mathrm{i}z\eta_1}=\exp(-z^2/2)$, which yields
$
\mathbb{L}_B(u)=\LL_{\tau} (u^2/2).
$
We can not see $\cL(B)\in T(\R_+)$ from this, and in fact $\cL(B)$
is in
$T(\R)$ but not in $T(\R_+)$ (see Bondesson \cite{Bondesson}, page 117).
Then using that $\eta$ is symmetric, we can apply (4) in Proposition
\ref{facts} and conclude that $\cL(V) \in T(\R)$ for suitable jump
heights of $\xi$.
\end{ex}

\begin{ex}[(The case when $\eta$ is a L\'evy subordinator and $\cL
(W_1)$ is a half normal
distribution)]
The $1/2$-stable subordinator $\eta$ and the standard half normal
random variable $\tau$ have densities, respectively, given by,
\[
f_{\eta_t}(x)=\frac{t}{2\sqrt{\uppi}} x^{-3/2} \mathrm{e}^{-t^2/2x}  \quad \mbox
{and} \quad   f_\tau(x)
=\sqrt{\frac{2}{\uppi}}\mathrm{e}^{-x^2/2},  \qquad x>0.
\]
These yield the density function of $\eta_\tau$ as
\[
f_{\eta_\tau}(x)=\int f_{\eta_y}(x)f_\tau(y)\,\mathrm{d}y =
\frac{1}{\sqrt{2}\uppi}\frac{x^{-1/2}}{1+x}.
\]
Interestingly, this is an $F$ distribution (see Sato \cite{Sa99}, page
46) and
since a random variable with an $F$ distribution is constricted to be the
quotient of two independent gamma random variables,
we have that $\cL(\eta_\tau)\in T(\R_+)$.
\end{ex}

\subsection*{Case 2: The process $\eta_t$ is a compound sum process}
Again we start with a condition for the convergence of
\eqref{integral}.
It can be shown similar to Proposition~\ref{convcondxirenewal}.

\begin{prop} \label{convcondetarenewal}
Suppose $\xi$ to be a L\'evy process and that
$\eta_t=\sum_{i=1}^{M_t}Y_i$ is a compound sum process with i.i.d. jump
heights $Y_i, i=1,2,\ldots,$ and i.i.d. waiting times $U_i,
i=1,2,\ldots.$ Then
\eqref{integral} converges a.s. to a finite random variable if and
only if
%
\begin{equation} \label{eq-convcondetarenewal}
\xi_t\to\infty \qquad\mbox{a.s.}  \quad    \mbox{and} \quad    \int_{(1,\infty)}
 \biggl( \frac{\log q}{A_\eta(\log q)}
 \biggr) P(|Y_1|\in \mathrm{d}q)<\infty
\end{equation}
for $A_\eta(x)=\int_{(0,x)} P(\xi_{U_1}>u)\,\mathrm{d}u.$
\end{prop}

In the same manner as before, we can now extend Proposition
\ref{CPPinThorin2} to the new setting and obtain the following result.

\begin{prop}\label{renewalinThorin2}
Let $\xi$ and $\eta$ be independent and assume
$\eta_t=\sum_{i=1}^{M_t}Y_i$ to be a compound renewal process with
i.i.d. jump
heights $Y_i, i=1,2,\ldots,$ and i.i.d. waiting times $U_i$ such that
\eqref{eq-convcondetarenewal} holds. Suppose that $\cL(Y_1)\in T(\R_+)$
and $\cL(\mathrm{e}^{-\xi_\tau})\in H(\R_+)$ for a random variable $\tau$
having the same
distribution as $U_1$
and being independent of $\xi$. Then
\[
\cL \biggl(\int_{(0,\infty)} \mathrm{e}^{-\xi_{t-}}\,\mathrm{d}\eta_t \biggr) \in T(\R_+).
\]
Furthermore, if $\cL(\mathrm{e}^{-\xi_\tau})\in\wt H(\R)$, $\cL(Y_1)\in
T(\R)$ and
$\cL(Y_1)$ is symmetric, then $\cL(V)\in T(\R)$.
\end{prop}

The following is a very simple example fulfilling the assumptions in
Proposition~\ref{renewalinThorin2}.

\begin{ex}[(The case when $\eta$ is a random walk and $\xi$ is a standard
Brownian motion with drift)]
Suppose $\xi_t=B_t+at$ is a standard Brownian motion with drift $a>0$
and $U_1$ is degenerated at 1.
Then $\cL(\mathrm{e}^{-\xi_\tau})=\cL(\mathrm{e}^{-a} \mathrm{e}^{-B_1})$ is a scaled log-normal
distribution and hence in $H(\R_+)$.
So for all GGC jump heights $\cL(Y_1)$, the exponential integral is GGC.
\end{ex}

\section*{Acknowledgements}
This research was carried out while Anita Behme was visiting Keio University
with financial support of ``Deutscher Akademischer Austauschdienst''. She
gratefully thanks for hospitality and funding. We thank the referee for a
careful reading of the manuscript and yielding suggestions.


\printhistory

\end{document}